\documentclass[11pt]{article}
\usepackage{amsmath,amssymb,amsfonts,mathtools,amsthm}
\usepackage{bm}
\usepackage{enumitem}
\usepackage{microtype}
\usepackage[a4paper,margin=1in]{geometry}
\usepackage[T1]{fontenc}
\usepackage{lmodern}
\theoremstyle{plain}
\newtheorem{theorem}{Theorem}[section]
\newtheorem{lemma}[theorem]{Lemma}
\newtheorem{proposition}[theorem]{Proposition}
\theoremstyle{definition}
\newtheorem{definition}[theorem]{Definition}
\theoremstyle{remark}

\newtheorem{example}[theorem]{Example}
\newcommand{\R}{\mathbb{R}}
\newcommand{\cS}{\mathcal S}
\DeclareMathOperator{\rank}{rank}
\DeclareMathOperator{\tr}{tr}

\DeclareMathOperator{\St}{St}

\title{The Sparse-Complement Conjecture for Two Distinct Eigenvalues}
\author{Binh T. Nguyen\thanks{Faculty of Mathematics and Computer Science, University of Science, Vietnam National University Ho Chi Minh City, Vietnam. Email: \texttt{ngtbinh@hcmus.edu.vn}.}}
\begin{document}
\maketitle

\begin{abstract}
For a graph $G$, let $q(G)$ be the minimum number of distinct eigenvalues
of a real symmetric matrix whose off-diagonal zero--nonzero pattern is $G$.
We prove that $q(G)=2$ whenever the complement $H=\overline G$ of a graph
of order $n\ge3$ has at most $n-3$ edges, resolving the sparse-complement
conjecture of Fallat and Mojallal. If $H$ contains a cycle, the proof
constructs a nontrivial orthogonal projector with the strong spectral
property and rank at most $|E(F)|-|V(F)|+4$, where $F$ is the union of
the cyclic components of $H$. The construction combines faithful
stress-free Gram representations, a one-coordinate lift for each excess
core edge, exact orthogonal-frame completion, and a Stiefel-manifold
pendant-extension argument.
\end{abstract}

\noindent\textbf{Keywords.}
inverse eigenvalue problem, graph complements, strong spectral property,
orthogonal projectors, Gram representations, Stiefel manifolds

\par\noindent\textbf{Mathematics Subject Classification (2020).}
05C50, 15A18, 15A29

\section{Introduction}\label{sec:intro}
The inverse eigenvalue problem of a graph asks which spectra occur among
real symmetric matrices whose off-diagonal zero--nonzero pattern is
prescribed by that graph.  The number of distinct eigenvalues is an
especially natural invariant because it depends on the permitted matrix
entries rather than on the adjacency matrix alone.  Let $q(G)$ denote its
minimum value for a graph $G$.  If $G$ is nonempty and admits a nontrivial
symmetric orthogonal projector, then $q(G)=2$; conversely any real symmetric
matrix with precisely two eigenvalues can be affinely normalized to such a
projector.

Fallat and Mojallal~\cite{FM23} posed the question whether removing at most
$n-3$ edges from a complete graph always preserves the possibility of a
realization with two distinct eigenvalues.  Barrett, Fallat, Furst,
Nasserasr, Rooney, and Tait~\cite{BF26} formulated the same statement as
Conjecture~1.1 and established, in particular, the case where the removed
edges form a bipartite graph.  The number $n-3$ is sharp: deleting $n-2$
edges incident with one vertex may produce a graph with a pendant vertex,
for which a two-eigenvalue realization is obstructed~\cite{BF26}.
Related results on symmetric orthogonal matrices of graph joins were
obtained by Levene, Oblak, and \v{S}migoc~\cite{LOS22,LOS25}; the strong
spectral property (SSP) was developed by Barrett and
coauthors~\cite{SSP17}.  The minimum vector rank of cycle complements,
including the special role of the four-cycle, was analyzed by Li,
Nathanson, and Phillips~\cite{LNP14}.

The construction works directly with the complement $H=\overline G$.
The cyclic components of $H$ are encoded by faithful Gram vectors, while its
tree components supply the rows needed to complete a Gram matrix to a
constant column-Gram (a tight frame).  The essential combinatorial identity
relates the number of tree components to the number of excess edges in the
cyclic core.  A single extra Gram coordinate enforces each excess edge,
without destroying the stress-free condition that later yields SSP.
The remaining forest vertices are then restored by a local transversality
argument on the Stiefel manifold.  Unlike a decomposition into balanced
joins, this mechanism does not require a prescribed partition of the
component orders.

\paragraph{Our contributions}
First, we resolve the sparse-complement conjecture of Fallat and
Mojallal. We show that every graph $G$ of
order $n\ge3$ with $|E(\overline G)|\le n-3$ has $q(G)=2$.
The result removes the bipartiteness hypothesis from the previously
known theorem at this edge threshold and extends the general bound
established in~\cite{BF26}. No restriction is imposed on the number
of cycles or the component orders of the complement.

The second contribution is a constructive spectral realization
mechanism with the strong spectral property. For complements with a
cyclic component, Theorem~\ref{thm:main} gives a nontrivial orthogonal
projector realizing $G$, with SSP and rank at most $k+4$, where
$k=|E(F)|-|V(F)|$ and $F$ is the union of the cyclic components of
$\overline G$. The construction combines faithful stress-free Gram
representations, one-coordinate lifts enforcing additional
orthogonality constraints, and exact orthogonal-frame completion.
The component-counting identity supplies the necessary completion
budget; a Stiefel-manifold transversality argument then restores
pendant trees while preserving the zero pattern and SSP.

These constructions build on established Gram-representation, orthogonal-
realization, and SSP techniques~\cite{BF26,SSP17,LOS22,LOS25}.
The new ingredient is their combination with a dimensional budget that also
covers the boundary case of $r-1$ available tree-component anchors.

\subsection{The main result and notation}
\label{subsec:main-result}

Let $G$ be a simple graph on $n$ vertices. We denote by
$\mathcal S(G)$ the set of real symmetric $n\times n$
matrices $A$ satisfying
\[
    A_{ij}\ne 0
    \quad\Longleftrightarrow\quad
    ij\in E(G), \qquad i\ne j.
\]
The minimum number of distinct eigenvalues of $G$ is
\[
    q(G)=\min_{A\in\mathcal S(G)}
          |\operatorname{spec}(A)|.
\]
Here, $\operatorname{spec}(A)$ denotes the set of
distinct eigenvalues of $A$.
A real symmetric matrix $P$ is an orthogonal projector
if $P^2=P$. When $0<\operatorname{rank}(P)<n$,
its spectrum consists precisely of $0$ and $1$.
Thus, constructing such a projector in $\mathcal S(G)$
establishes $q(G)\le 2$.

We write $H=\overline G$ for the complement of $G$.
For a matrix $P\in\mathcal S(G)$, the edges of $H$
are precisely the off-diagonal zero positions of $P$.
Throughout, isolated vertices are regarded as
tree components.
We establish the sparse-complement conjecture
together with a stronger realization statement
for complements containing cycles.

\begin{theorem}[Sparse-complement conjecture]
\label{thm:main}
Let $G$ be a simple graph of order $n\ge 3$,
and let $H=\overline G$. If
\[
    |E(H)|\le n-3,
\]
then $q(G)=2$.

Moreover, suppose that $H$ contains a cycle,
and let $F$ be the union of all connected components
of $H$ containing at least one cycle. Set
\[
    k=|E(F)|-|V(F)|.
\]
Then there exists an orthogonal projector
$P\in\mathcal S(G)$ with the strong spectral
property and
\[
    3\le \operatorname{rank}(P)\le k+4.
\]
\end{theorem}

The proof uses faithful Gram representations
and a nondegeneracy condition that ensures
the strong spectral property after frame
completion.

\begin{definition}[Faithful Gram representations]
\label{def:stress}
Let $F$ be a graph on $s$ vertices, and let
$M\in\mathbb R^{s\times r}$ have row vectors
$m_1^\top,\ldots,m_s^\top$.

We call $M$ a \emph{faithful Gram representation}
of $F$ if every row is nonzero and
\[
    m_i^\top m_j=0
    \quad\Longleftrightarrow\quad
    ij\in E(F), \qquad i\ne j.
\]
The representation is \emph{nonparallel} if
no two distinct row vectors are collinear.

We call $M$ \emph{stress-free} if
\begin{equation}\label{eq:stress}
\begin{aligned}
    X&=X^\top, \qquad X_{ii}=0,\\
    X_{ij}&=0 \quad\text{if }ij\notin E(F),\\
    XM&=0
\end{aligned}
\qquad\Longrightarrow\qquad X=0.
\end{equation}
\end{definition}

\begin{definition}[Strong spectral property]
\label{def:ssp}
A real symmetric matrix $P$ has the
\emph{strong spectral property} (SSP) if
the only real symmetric matrix $X$ satisfying
\[
    X_{ii}=0,\qquad
    X_{ij}=0 \text{ whenever }P_{ij}\ne0,
    \qquad XP=PX
\]
is $X=0$.
\end{definition}

The proof proceeds in three stages. We first
construct faithful, nonparallel, stress-free
Gram representations of cyclic components
(Sections~\ref{sec:core}--\ref{sec:generic}).
We then complete these representations to
orthogonal projectors with SSP
(Section~\ref{sec:completion}).
Finally, we recover the remaining tree vertices
through successive pendant extensions on a
Stiefel manifold (Section~\ref{sec:pendant}).
The main theorem is proved in
Section~\ref{sec:proof}.

\section{Related work and positioning}
\label{sec:related}

The minimum number of distinct eigenvalues of a graph has
been studied from both extremal and structural perspectives.
Barrett et al.~\cite{BFFK23} determined the minimum number
of edges in a connected graph admitting a two-eigenvalue
realization, while related classifications address
low-degree regular graphs~\cite{BFFR23}. In the complementary
direction, Fallat and Mojallal~\cite{FM23} conjectured that
$q(G)=2$ whenever $|E(\overline G)|\le n-3$.
Barrett et al.~\cite{BF26} established this conclusion
for bipartite complements and proved the general bound
$q(G)=2$ when
$|E(\overline G)|\le \lfloor n/2\rfloor-1$.
Their work also employs Gram representations and the
strong spectral property in constructions for certain
cyclic complements.

Several techniques for constructing symmetric matrices
with few distinct eigenvalues are closely related to
our approach. Levene, Oblak, and \v{S}migoc
developed orthogonal realizations for graph joins
\cite{LOS22,LOS25} and investigated generic eigenvector
constructions~\cite{LOS24}. Abiad et al.~\cite{AFK23}
studied symmetric-matrix bordering, while Culver and
Kempton~\cite{CK25} obtained further two-eigenvalue
realizations through graph products. The strong spectral
property provides a transversality framework for
prescribed matrix patterns~\cite{SSP17}.
These results supply important structural and
methodological foundations, although their stated
hypotheses do not directly cover arbitrary complements
with at most $n-3$ edges.

Related investigations of restricted graph classes
include threshold graphs~\cite{FM22Threshold,
Allem25Threshold,BFP26Threshold} and
cographs~\cite{AFH26Cograph}, for which uniform
upper bounds on $q(G)$ have been obtained.
These results further illustrate the influence of
graph structure on attainable eigenvalue
multiplicities.

Li, Nathanson, and Phillips~\cite[Section~4]{LNP14}
studied orthogonal representations of complements of
cycles and unicyclic graphs, together with extension
methods for pendant vertices and ears and results for
structured multicyclic families such as necklaces.
Those constructions concern low-dimensional vector
representations; they do not impose the simultaneous
requirements of a fixed-rank orthogonal projector and
SSP. Building on the cycle representations, our proof
combines Gram-coordinate lifting, orthogonal-frame
completion, and Stiefel transversality to handle
arbitrary cyclic components of sparse complements,
including the boundary case of only $r-1$ available
tree-component anchors. The resulting realization
has SSP and an explicit rank bound.

\section{Core, budget, and Gram representation}\label{sec:core}
\begin{lemma}[Tree budget]\label{lem:budget}
Let $F$ be the disjoint union of all connected components of $H$ containing a cycle. Let $T$ count the tree components of $H$, including isolates, and put $k=e(F)-v(F)$. Then,
\[e(H)=v(H)+k-T,\qquad e(H)\le v(H)-3\Longrightarrow T\ge k+3.\]
\end{lemma}
\begin{proof}
Every cyclic component $C$ contributes $v(C)+(e(C)-v(C))$ edges and each tree of order $t$ contributes $t-1$. Sum over components.
\end{proof}

\begin{lemma}[Unicyclic cores]\label{lem:seed}
Every disjoint union of unicyclic graphs containing at least one cycle has a faithful, nonparallel, stress-free Gram representation of column rank $d_0\in\{3,4\}$. Dimension $3$ is sufficient if the chosen base cycles have no $C_4$; dimension $4$ is sufficient in all cases.
\end{lemma}
\begin{proof}
For $\ell\ge5$, Li, Nathanson, and
Phillips~\cite[Proposition~4.2 and the ensuing discussion]{LNP14}
provide an orthogonal representation of $\overline{C_\ell}$
in $\R^3$ by nonzero, pairwise linearly independent vectors.
Their convention places orthogonality at nonedges of the
represented graph. Thus, for the same vectors,
\[
 m_i^\top m_j=0
 \quad\Longleftrightarrow\quad
 ij\in E(C_\ell)\qquad(i\ne j),
\]
which is precisely a faithful, pairwise nonparallel Gram
representation of $C_\ell$ under our convention. For $C_3$,
the standard basis of $\R^3$ gives such a representation.
For $C_4$, take
\[
 (1,0,0,0),\quad (0,0,1,0),\quad
 (1,1,0,0),\quad (0,0,1,1).
\]
These vectors form a faithful, pairwise nonparallel
representation spanning $\R^4$, consistent with the
four-cycle distinction in~\cite[Section~4.1]{LNP14}.
Work in $\R^d$, where $d=4$ if a component contains a
four-cycle and $d=3$ otherwise, embedding
three-dimensional seeds when needed.

Each remaining vertex of a unicyclic component can be added in
parent-before-child order from the base cycle. Suppose that the vectors
of all earlier vertices are nonzero and pairwise nonparallel, and that
the new vertex $z$ has parent $p$. Choose $m_z\in m_p^\perp$.
For every earlier $i\ne p$, the additional condition
$m_z^\top m_i=0$ excludes the proper linear subspace
$m_p^\perp\cap m_i^\perp$ of $m_p^\perp$, because $m_i$ is not
parallel to $m_p$. The forbidden zero vector and the finitely many
collinear directions also form proper subsets of $m_p^\perp$, which
has dimension $d-1\ge2$. Their union does not exhaust this space.
Thus $m_z$ can be chosen nonzero and nonparallel to all earlier rows,
orthogonal precisely to its parent.

For distinct components, independently rotate their representations
in $O(d)$. For fixed nonzero vectors $a,b\in\R^d$, neither the
equality $a^\top Rb=0$ nor the condition $Rb\in\R a$ holds
identically on either connected component of $O(d)$: rotations can
move $Rb$ through the unit sphere, and $d\ge3$. Their exceptional
sets are proper real-analytic subsets. Since there are finitely many
cross-component pairs, the rotations can be chosen so that all
cross-component inner products are nonzero and all rows remain
pairwise nonparallel. At least one base cycle spans $\R^3$ (or
$\R^4$ in the presence of a four-cycle), so the resulting row-span
has dimension $d$.

Finally, suppose $X$ is symmetric, supported on the edges of the
resulting disjoint union, has zero diagonal, and satisfies $XM=0$.
For a leaf $z$ with parent $p$, row $z$ of this equation is
$X_{zp}m_p=0$, whence $X_{zp}=0$. Remove all such leaf coefficients
recursively. At a cycle vertex, the two remaining neighbor vectors
are nonparallel and hence linearly independent. The corresponding
row equation forces both edge coefficients to vanish. This applies
componentwise and gives $X=0$.
\end{proof}

\begin{lemma}[Exact one-edge lift]\label{lem:lift}
Suppose $M$ is faithful, nonparallel, and stress-free for $F$ and $uv$ is a nonedge of $F$. Append the column $z=e_u-c e_v$, $c=m_u^\top m_v\ne0$. Then, $[M\ z]$ is faithful, nonparallel, and stress-free for $F+uv$.
\end{lemma}
\begin{proof}
The rank-one Gram perturbation $zz^\top$ changes only the off-diagonal $uv$-entry, from $c$ to $c-c=0$. Projection onto the old coordinates preserves nonparallelity. If $X[M\ z]=0$ and $X$ is symmetric and supported on $E(F+uv)$, then $(Xz)_u=-cX_{uv}=0$, hence, $X_{uv}=0$. Now $XM=0$ with support in $E(F)$, so $X=0$.
\end{proof}

\begin{lemma}[Effective-rank core]\label{lem:core}
For the cyclic core $F$ with $k=e(F)-v(F)\ge0$, there exists a faithful, nonparallel, stress-free, full-column-rank $M\in\R^{s\times r}$ such that $3\le r\le k+4$.
\end{lemma}
\begin{proof}
In each cyclic component, retain a spanning tree together with one
additional edge. The resulting spanning unicyclic subgraphs admit a
joint representation of ambient dimension $d_0\le4$ by
Lemma~\ref{lem:seed}. Each component contributes exactly
$|E(C)|-|V(C)|$ omitted edges, giving $k$ omitted edges in total.
Add them one at a time using Lemma~\ref{lem:lift}, and denote the
resulting matrix by $B\in\R^{s\times(d_0+k)}$.

Let $r=\rank(B)$ and let $R\in\R^{(d_0+k)\times r}$ have orthonormal
columns spanning the row space of $B$. Then
\[
 B=(BR)R^\top,\qquad BB^\top=(BR)(BR)^\top.
\]
Consequently $M=BR$ has the same Gram matrix and hence the same
faithful pattern and nonparallel rows. Moreover, for every matrix
$X$ of compatible size, $XB=0$ if and only if $X(BR)=0$;
stress-freeness is therefore preserved. The cycle seed already spans
at least three dimensions, while adding $k$ columns increases rank
by at most $k$. Thus $3\le r\le d_0+k\le k+4$, as required.
\end{proof}

\begin{lemma}[Positive row scaling preserves the properties]\label{lem:scalepreserve}
If $M$ is faithful, nonparallel, and stress-free, so is $DM$ for every positive diagonal $D$.
\end{lemma}
\begin{proof}
Zero/nonzero inner products are scaled by positive factors and nonparallelity persists. If $XDM=0$, then $(DXD)M=0$; the matrix $DXD$ is symmetric with the same permitted support. Stress-freeness implies $DXD=0$, hence, $X=0$.
\end{proof}

\section{Genericity lemmas}\label{sec:generic}
\begin{lemma}[Simple spectrum and avoidance by positive scaling]\label{lem:scaling}
Let $m_1,\dots,m_s$ be nonzero, pairwise nonparallel rows spanning $\R^r$, $r\ge3$. Put $Q(w)=\sum_i w_i m_i m_i^\top$ with $w_i>0$.
\begin{enumerate}[label=(\roman*),leftmargin=*]
\item The weights for which $Q(w)$ has distinct positive eigenvalues form a dense open subset of the positive orthant.
\item If some $m_u,m_v$ are nonorthogonal (and hence nonparallel), there exist positive weights giving a simple largest eigenvalue whose unit eigenvector is parallel to none of the $m_i$. The two conditions can hold simultaneously.
\end{enumerate}
\end{lemma}
\begin{proof}
(i) Select linearly independent vectors $m_{i_1},\ldots,m_{i_r}$.
For $j=1,\ldots,r$, the positive semidefinite matrix
$A_j=\sum_{\ell=1}^j t_\ell m_{i_\ell}m_{i_\ell}^\top$ has rank $j$.
Inductively choose $t_j>0$ sufficiently small that the $j-1$
previously positive, distinct eigenvalues remain distinct and separated
from zero.  The additional positive eigenvalue is then smaller than
all of them.  At $j=r$ the spectrum is positive and simple.
Adding arbitrarily small positive weights to the remaining rows preserves
simplicity.  Consequently, the discriminant of the characteristic
polynomial of $Q(w)$ is a nonzero polynomial in the weights.  Its
nonvanishing locus is open and dense in $(0,\infty)^s$.

(ii) Write $a=m_u$, $b=m_v$, so $a^\top b\ne0$ and $a,b$
are linearly independent, and put $A(t)=taa^\top+bb^\top$.
For $t>0$ the two positive eigenvalues of $A(t)$ are distinct:
in the basis of the plane $\operatorname{span}\{a,b\}$ obtained by taking
the first basis vector along $a$, the off-diagonal entry is nonzero.
Let $h(t)$ span its top eigenspace.  Neither $a^\top h(t)$ nor
$b^\top h(t)$ vanishes.  Indeed, if $a^\top h(t)=0$, the
eigenvalue equation forces $h(t)$ to be parallel to $b$, a
contradiction with $a^\top b\ne0$; the other case is analogous.
If a fixed line $\R h$ were the top eigenspace at two distinct
parameters $t_1,t_2$, subtracting the eigenvalue equations gives
$
 (t_1-t_2)a(a^\top h)
   =(\lambda_1(t_1)-\lambda_1(t_2))h.
$

The left-hand coefficient $a^\top h$ is nonzero, so $h$ would be
parallel to $a$, contrary to $b^\top h\ne0$ and the eigenvalue
equation for $A(t)$.  Thus, any prescribed line occurs for at most
one $t$.  Avoiding the finitely many lines $\R m_i$ leaves a value
$t>0$ for which the simple top eigendirection of $A(t)$ is
admissible.  Sufficiently small positive weights on all remaining
rows preserve both properties.  This yields a nonempty open set of
admissible positive weights; by (i) it intersects the dense open
simple-spectrum locus.
\end{proof}

\begin{lemma}[Generic orthogonal mixing]\label{lem:mix}
Let $W_0\in\R^{T\times r}$ have rank at least $2$, and suppose $C=W_0W_0^\top$ is not scalar. Let $m_1,\dots,m_s$ be nonzero rows such that $W_0m_i\ne0$ for all $i$. Then, there exists $R\in O(T)$ such that all rows of $W=RW_0$ are nonzero and pairwise nonparallel, no row of $W$ is parallel to any $m_i$, and every $w_a^\top w_b$ ($a\ne b$) and $w_a^\top m_i$ is nonzero.
\end{lemma}
\begin{proof}
Write $L=W_0^\top:\R^T\to\R^r$ and $C=W_0W_0^\top$.
For $R\in O(T)$, row $a$ of $RW_0$ equals
$L(R^\top e_a)^\top$.  Each forbidden equality is the vanishing
of a polynomial (hence analytic) function of the entries of $R$:
$\|w_a\|^2$,$
 \|w_a\wedge w_b\|^2$, 
 $\|w_a\wedge m_i\|^2$,
 $w_a^\top w_b$, and $w_a^\top m_i$.

We verify that none of these functions vanishes identically on $O(T)$.
The first is nontrivial because $L\ne0$.  Since $\rank L\ge2$,
there are orthonormal $x,y\in\R^T$ for which $Lx,Ly$ are
linearly independent: take a unit $x$ with $Lx\ne0$, then choose
$y\perp x$ with $Ly\notin\R Lx$; otherwise $\operatorname{im}L$
would be one-dimensional.  This proves nontriviality of the second.
For fixed $m_i\ne0$, $\operatorname{im}L$ is not contained in
$\R m_i$, which proves nontriviality of the third.  Since $C$ is
symmetric and nonscalar, it has orthogonal eigenvectors with distinct
eigenvalues; rotating their span by $\pi/4$ gives orthonormal
$x,y$ with $x^\top Cy\ne0$, proving the fourth.  Finally,
$W_0m_i\ne0$ implies that $x^\top W_0m_i\ne0$ for some unit
$x$, proving the fifth.  In every case the chosen one or two domain
vectors can be extended to an orthonormal basis, so they can be
realized as specified rows of $R$.

Here $\|u\wedge v\|^2=\|u\|^2\|v\|^2-(u^\top v)^2$,
so each nonparallelity constraint is polynomial. The zero set of a
nonzero real-analytic function has empty interior on each connected
component on which it is nonzero.  A diagonal sign
matrix with exactly one negative entry exchanges the two components
of $O(T)$; it multiplies any listed expression by $\pm1$ or leaves it
unchanged.  Thus, nontriviality holds on both components.  Each
exceptional set is therefore nowhere dense in $O(T)$.  Their finite
union is not all of $O(T)$, which proves the assertion.
\end{proof}

\section{Exact frame completion and SSP}\label{sec:completion}
\begin{lemma}[Universal-vertex completion]\label{lem:complete}
Let $M\in\R^{s\times r}$ be full column rank, faithful, nonparallel and stress-free, with $r\ge3$. If $T\ge r$, there is a rank-$r$ SSP projector in $\cS(\overline{F\sqcup TK_1})$ with nonzero pairwise nonparallel Gram rows. If $T=r-1$ and $F$ is not complete, the same holds.
\end{lemma}
\begin{proof}
By Lemmas~\ref{lem:scalepreserve}--\ref{lem:scaling}, positively rescale the rows to make $Q=M^\top M$ have simple positive spectrum. If $T\ge r$, take $\alpha>\lambda_1(Q)$ and $S=\alpha I_r-Q\succ0$. Factor $S=W_0^\top W_0$ with $W_0\in\R^{T\times r}$ of rank $r$. If $T=r-1$, a noncomplete $F$ has a nonedge, hence two nonorthogonal nonparallel Gram rows; apply Lemma~\ref{lem:scaling}(ii), choose $\alpha=\lambda_1(Q)$, and factor $S=\alpha I_r-Q\succeq0$ with $W_0\in\R^{(r-1)\times r}$ of rank $r-1$. Its kernel is $\R h$, with $h$ the top eigenvector of $Q$; no $m_i$ lies in this kernel. In both cases $W_0m_i\ne0$.

The row Gram matrix $W_0W_0^\top$ is nonscalar: in the full-rank case this follows from $\rank(W_0)=r$ when $T>r$, and from the nonscalar simple-spectrum $S$ when $T=r$; in the tight case its positive eigenvalues are the distinct numbers $\alpha-\lambda_2,\ldots,\alpha-\lambda_r$. Lemma~\ref{lem:mix} supplies $W=RW_0$ with the required zero/nonzero and nonparallel pattern. Let
\[ U=\alpha^{-1/2}\begin{bmatrix}M\\W\end{bmatrix},\qquad P=UU^\top.\]
Since $U^\top U=I_r$, $P$ is a rank-$r$ projector. Its off-diagonal zeros are exactly the edges of $F$.

Let $X$ be a symmetric SSP certificate for $P$. The zero pattern implies $X=\operatorname{diag}(X_F,0_T)$ with $X_F$ supported on $E(F)$ and zero diagonal. From $XP=PX$, we obtain $XU=UB$ for symmetric $B=U^\top XU$; hence, $X_FM=MB$ and $WB=0$. In the full-rank case, $B=0$. In the tight case, $\ker W=\R h$.  As $B$ is symmetric,
$WB=0$ implies $\operatorname{im} B\subseteq\R h$, hence
$B=\tau hh^\top$ for a unit vector $h$.  Moreover,
\[
 \tr(X_FMM^\top)=\sum_{i,j}(X_F)_{ij}\,m_i^\top m_j=0,
\]
where the diagonal terms vanish, and every off-diagonal term has
$(X_F)_{ij}\ne0$ only on an edge of $F$, where the Gram
product vanishes.  Therefore,
\[0=\tr(X_FMM^\top)=\tr(M^\top M B)=\alpha\tau.\]
Thus, $B=0$ in either case, and stress-freeness gives $X_F=0$. Hence, $P$ has SSP.
\end{proof}

\begin{lemma}[The tight case is noncomplete]\label{lem:clique}
Suppose $F$ is a nonempty union of cyclic components, $M$ is full column rank and faithful in dimension $r$, and $T\ge k+3$ with $k=e(F)-v(F)$. If $T=r-1$, then $F$ cannot be complete.
\end{lemma}
\begin{proof}
A complete $F$ must be a single clique $K_s$, $s\ge3$. Faithfulness forces its $s$ Gram rows mutually orthogonal and nonzero, so $r=s$. But
\[T\ge\binom{s}{2}-s+3=s+\frac{(s-2)(s-3)}2\ge s=r,\]
contradicting $T=r-1$.
\end{proof}

\section{A Stiefel pendant-extension theorem}\label{sec:pendant}
\begin{lemma}[Exact pendant extension preserving SSP]\label{lem:pendant}
Let $P=UU^\top$ be an SSP projector in $\cS(\overline H)$, with $U\in\R^{n\times r}$, $U^\top U=I_r$, $r\ge3$, and nonzero pairwise nonparallel rows. Form $H^+$ by adjoining a new vertex $z$ adjacent only to $p$ in $H$. Then $\overline{H^+}$ has a rank-$r$ SSP projector with nonzero pairwise nonparallel frame rows.
\end{lemma}
\begin{proof}
Use the convention
$\St(r,N)=\{V\in\R^{N\times r}:V^\top V=I_r\}$.
Put $V_0=[U;0]\in\St(r,n+1)$ and let
\[
 f:\St(r,n+1)\longrightarrow\R^{E(H^+)},\qquad
 f(V)=(v_i^\top v_j)_{ij\in E(H^+)}.
\]
Then, $f(V_0)=0$.  For a tangent matrix $Y=(y_i)_{i\in V(H^+)}$,
\[
 (Df(V_0)Y)_{ij}=y_i^\top (V_0)_j+(V_0)_i^\top y_j,
 \qquad V_0^\top Y+Y^\top V_0=0.
\]
Suppose a linear combination of these coordinate differentials
vanishes on the tangent space.  Encode its coefficients by a symmetric
matrix $X$ with zero diagonal and support in $E(H^+)$.  The ambient
gradient is $2XV_0$, so orthogonality to the tangent space gives
$XV_0=V_0B$ for some symmetric $B$: the normal space at $V_0$
is exactly $\{V_0B:B=B^\top\}$.
The last row reads $x_{zp}u_p=0$ and forces $x_{zp}=0$ because
$u_p\ne0$.  Writing $X_0$ for the old principal block, we have
$X_0U=UB$, whence $X_0$ commutes with $UU^\top$.  Since $X_0$
has zero diagonal and is supported on $E(H)$, SSP gives $X_0=0$.
Thus, the differentials are independent, and $Df(V_0)$ is onto.

Choose $w\in u_p^\perp$ such that
$w^\top u_i\ne0$ for every $i\ne p$, and $w$ is not parallel
to any $u_i$.  Such a vector exists because $\dim u_p^\perp=r-1\ge2$,
none of the old rows is parallel to $u_p$, and each forbidden
orthogonality or direction is a proper linear subset of
$u_p^\perp$.  The matrix $Y=[0;w]$ is tangent at $V_0$ and
$Df(V_0)Y=0$.  Since $Df(V_0)$ is onto, the regular-level-set theorem shows that
$f^{-1}(0)$ is a smooth submanifold near $V_0$, with tangent space
$\ker Df(V_0)$. Every tangent vector to a smooth manifold is the
velocity of a smooth curve on the manifold. Hence there is a smooth
curve $V(t)\in f^{-1}(0)$ with $V(0)=V_0$ and $V'(0)=Y$.
In particular, its new row is $tw+o(t)$, so for every sufficiently
small $t>0$ it has precisely the required nonzero inner products
with old rows other than $p$; orthogonality to $p$ holds identically
because $f(V(t))=0$.  Old nonzero products remain nonzero, and
nonzero/nonparallel row conditions persist for all sufficiently
small $t>0$ (for the new row, divide by $t$ and pass to $w$).

Finally, on the fixed vector space of symmetric zero-diagonal
matrices supported on $E(H^+)$, the linear map
$X\mapsto X(VV^\top)-(VV^\top)X$ is injective at $V_0$ by the
same normal-space calculation above.  Its matrix entries depend
continuously on $V$, and full column rank is open.  It remains
injective at $V(t)$ for small $t>0$, which is precisely SSP.
\end{proof}

\section{Proof of the main theorem}\label{sec:proof}

\begin{proof}[Proof of Theorem~\ref{thm:main}]
Put $H=\overline G$ and $n=|V(H)|$. We first distinguish
whether $H$ contains a cycle.

\smallskip\noindent
\emph{Case 1: $H$ is a forest.}
A forest is bipartite. Thus, the hypothesis $|E(H)|\le n-3$
allows us to apply the bipartite-complement theorem
\cite[Theorem~3.7]{BF26}, which gives $q(G)=2$. The additional
projector and SSP assertion in the statement of the present
theorem concerns only complements containing a cycle, so
nothing further is required in this case.

\smallskip\noindent
\emph{Case 2: $H$ has a cyclic component.}
Let $F$ be the disjoint union of all connected components
of $H$ that contain a cycle, and write $s=|V(F)|$.
Let $T$ be the number of remaining (tree) components,
including isolated vertices, and set
\[
    k=|E(F)|-|V(F)|.
\]
Every component of $F$ has at least as many edges as
vertices; hence $k\ge0$. By Lemma~\ref{lem:budget}
and the hypothesis $|E(H)|\le n-3$,
\begin{equation}\label{eq:main-budget}
    T\ge k+3\ge3.
\end{equation}
Lemma~\ref{lem:core} supplies a matrix
$M\in\R^{s\times r}$ whose rows form a faithful,
pairwise nonparallel, stress-free Gram representation
of $F$, with full column rank and
\begin{equation}\label{eq:main-rank}
    3\le r\le k+4.
\end{equation}
Combining \eqref{eq:main-budget} and
\eqref{eq:main-rank} gives
\begin{equation}\label{eq:main-anchors}
    T\ge k+3\ge r-1.
\end{equation}
This is precisely the number of anchors needed
for the frame completion below.
Choose one root in each of the $T$ tree components
of $H$. Delete all other vertices in those components,
leaving the induced graph
\[
    H_0=F\,\dot\cup\,TK_1.
\]
We next construct an orthogonal projector with graph
$\overline{H_0}$.
If $T\ge r$, apply the first assertion of
Lemma~\ref{lem:complete} to $M$. Otherwise,
\eqref{eq:main-anchors} forces $T=r-1$.

In that boundary case, Lemma~\ref{lem:clique}
shows that $F$ is not complete, so the second
assertion of Lemma~\ref{lem:complete} applies.
Thus, in either case, there is a matrix
\[
    U_0\in\R^{(s+T)\times r},
    \qquad U_0^\top U_0=I_r,
    \qquad P_0=U_0U_0^\top,
\]
for which $P_0$ has rank $r$, belongs to
$\mathcal S(\overline{H_0})$, and has SSP.
Moreover, all rows of $U_0$ are nonzero and
pairwise nonparallel.
In particular, for distinct vertices $i,j$ of $H_0$,
\[
    (P_0)_{ij}=0
    \quad\Longleftrightarrow\quad
    ij\in E(H_0);
\]
the equivalence, not merely the forward implication,
is provided by the faithful frame completion.

It remains to recover the deleted vertices.
Root each original tree component of $H$ at its
retained vertex and order its nonroot vertices
so that every vertex follows its parent.
Starting from $H_0$, add these vertices one at a time
in that order. At each step, the new vertex has
exactly one neighbor among the vertices already
present, namely its parent: all other edges incident
to it in its original tree lead to descendants
not yet added. No edge joins different components
of $H$.
Consequently, every step is an instance of the
pendant extension in Lemma~\ref{lem:pendant}.
Inductively, that lemma produces a rank-$r$
orthogonal projector with SSP and with nonzero,
pairwise nonparallel frame rows for the complement
of the graph constructed so far. Its full
zero--nonzero pattern is preserved at each extension.

After finitely many steps the graph is $H$,
so we obtain
\[
    P=P^\top=P^2\in\mathcal S(\overline H)
      =\mathcal S(G),
    \qquad \rank(P)=r,
\]
and $P$ has SSP. Equation~\eqref{eq:main-rank}
gives the claimed rank bound
\[
    3\le\rank(P)\le k+4.
\]
Finally, $M$ has full column rank $r$, hence $s\ge r$.
Also $T\ge3$ by \eqref{eq:main-budget}, so
\[
    0<r\le s<s+T\le n.
\]
Therefore, $P$ is neither the zero nor the identity
projector, and its spectrum contains both $0$
and $1$ and no other values. This proves $q(G)\le2$.

Conversely, the inequality $|E(H)|\le n-3$
and $n\ge3$ imply
\[
    |E(H)|<\binom n2,
\]
so $G$ has at least one edge.
A real symmetric matrix with only one distinct
eigenvalue is scalar, and hence has no nonzero
off-diagonal entries. No such matrix can belong
to $\mathcal S(G)$; thus $q(G)\ge2$.

We conclude that
\[
    {q(G)=2},
\]
as required.
\end{proof}

\section{Sharpness and an exact rational realization}
\label{sec:sharpness}

We first show that the edge threshold in
Theorem~\ref{thm:main} is best possible.

\begin{proposition}\label{prop:sharpness}
For every $n\ge 4$, there exists a graph $G$
of order $n$ such that
\[
    |E(\overline G)|=n-2
    \qquad\text{and}\qquad
    q(G)>2.
\]
Consequently, the bound $n-3$ in
Theorem~\ref{thm:main} cannot be increased.
\end{proposition}

\begin{proof}
Let
$H=K_{1,n-2}\,\dot\cup\,K_1$ and $G=\overline H$. 
Clearly, $|V(H)|=n$ and $|E(H)|=n-2$.
Let $u$ be the center of the star, $v$ its
isolated vertex, and $w$ any leaf of the star.
In $G$, the vertex $u$ is adjacent only to $v$,
whereas $v$ is adjacent to every other vertex.
Thus, $u,v,w$ form a unique shortest path
of length two between $u$ and $w$.
The unique-shortest-path bound
\cite[Observation~1.3]{BF26} gives $q(G)\ge3$.
Hence, the threshold $n-3$ is sharp.
\end{proof}

The projector-rank bound obtained in
Theorem~\ref{thm:main} is a bound for the
present construction; its optimality is
not required for the sharpness assertion.
We next give an exact rational realization
for a complement containing a component with
more than one independent cycle.

\begin{example}\label{ex:rational}
Let
\[
    H=K_4\,\dot\cup\,5K_1.
\]
Then, $n=9$ and
\[
    |E(H)|=6=n-3.
\]
In particular, $H$ is not a pseudoforest.
Set
\[
    D=\operatorname{diag}(1,4,7,8),
    \qquad
    E=
    \begin{bmatrix}
      8&0&0&0\\
      0&7&0&0\\
      0&0&4&0\\
      0&0&0&1\\
      0&0&0&0
    \end{bmatrix}.
\]
For $j=1,2,3$, define
\[
    H_j=I_5-
    \frac{2v^{(j)}(v^{(j)})^\top}
         {(v^{(j)})^\top v^{(j)}},
\]
where
\[
\begin{aligned}
    v^{(1)}&=(7,7,1,5,8)^\top,\\
    v^{(2)}&=(7,5,8,6,4)^\top,\\
    v^{(3)}&=(3,5,3,2,5)^\top.
\end{aligned}
\]
Put
\[
    R=H_3H_2H_1,\qquad
    A=\begin{bmatrix}D\\RE\end{bmatrix},
    \qquad
    P=\frac1{65}AA^\top.
\]
Each $H_j$ is a rational orthogonal matrix,
and hence $R^\top R=I_5$. Moreover,
\[
    A^\top A
    =D^2+E^\top E
    =65I_4.
\]
Consequently,
\[
    P^\top=P,\qquad P^2=P,\qquad
    \operatorname{rank}(P)=4.
\]
Direct evaluation over $\mathbb Q$ gives
\[
    P_{ij}=0
    \quad\Longleftrightarrow\quad
    1\le i<j\le4
    \qquad (i<j).
\]
Thus, $P\in\mathcal S(\overline H)$.
It remains to verify the strong spectral
property. Write
\[
    P=\begin{bmatrix}
        P_{11}&P_{12}\\
        P_{12}^\top&P_{22}
      \end{bmatrix},
    \qquad
    P_{12}=\frac1{65}DE^\top R^\top.
\]
Since $D$ is invertible, $E$ has column
rank four, and $R$ is orthogonal,
\[
    \operatorname{rank}(P_{12})=4.
\]
Suppose that $X$ is symmetric, has zero
diagonal, vanishes on the edges of
$\overline H$, and satisfies $XP=PX$.
The prescribed support of $X$ forces
\[
    X=\begin{bmatrix}
        X_0&0\\
        0&0
      \end{bmatrix},
    \qquad X_0\in\mathbb S^4.
\]
The upper-right block of $XP=PX$ yields
\[
    X_0P_{12}=0.
\]
Since $P_{12}$ has full row rank,
$X_0=0$, and therefore $X=0$.
This proves that $P$ has SSP.
Hence, $\overline H$ admits a rational
rank-four orthogonal projector with
the strong spectral property.
\end{example}

\section{Conclusion and open questions}
\label{sec:conclusion}

We have established that every graph $G$ of order $n\ge 3$
whose complement has at most $n-3$ edges satisfies $q(G)=2$,
resolving the sparse-complement conjecture. The proof is
constructive in structure: a component-counting identity
provides the necessary dimensional budget, Gram-coordinate
lifting handles excess edges of cyclic components, and
orthogonal-frame completion produces projectors with the
strong spectral property. A transversality argument on
the Stiefel manifold then extends the construction through
successive pendant attachments.

The method also yields a rank bound for the resulting
projectors when the complement contains a cycle. Writing
$F$ for the union of its cyclic components and
$k=|E(F)|-|V(F)|$, the construction gives
\[
    3\le \operatorname{rank}(P)\le k+4.
\]
This bound arises from the Gram representation and
completion procedure and is not asserted to be optimal.

Two questions naturally arise. First, what is the minimum
rank of a projector realizing a prescribed sparse-complement
pattern, possibly under the additional requirement of the
strong spectral property? Second, can the existential
construction be made algorithmic, with effective complexity
bounds or exact-arithmetic guarantees? These questions
concern the quantitative and computational aspects of
two-eigenvalue realizations beyond their existence.

\section*{Declarations}
\noindent\textbf{Use of artificial intelligence.}
Generative AI tools assisted with literature
exploration, proof checking, and manuscript writing and editing. The authors assume
responsibility for all content.

\bibliographystyle{plain}
\bibliography{references}

@article{FM23,
 author = {Shaun Fallat and Seyed Ahmad Mojallal},
 title = {Spectral Applications of Vertex-Clique Incidence Matrices Associated with a Graph},
 journal = {Mathematics},
 year = {2023}, volume = {11}, number = {16}, pages = {3595},
 doi = {10.3390/math11163595}
}

@article{BF26,
 author = {Wayne Barrett and Shaun Fallat and Veronika Furst and Shahla Nasserasr and Brendan Rooney and Michael Tait},
 title = {Graphs with Bipartite Complement that Admit Two Distinct Eigenvalues},
 journal = {Electronic Journal of Linear Algebra}, year = {2026}, volume = {42}, pages = {146--161},
 doi = {10.13001/ela.2026.9443}
}

@article{LNP14,
 author = {Xiaowei Li and Michael Nathanson and Rachel Phillips},
 title = {Minimum Vector Rank and Complement Critical Graphs},
 journal = {Electronic Journal of Linear Algebra}, year = {2014}, volume = {27}, pages = {100--123},
 doi = {10.13001/1081-3810.1606}
}

@article{SSP17,
 author = {Wayne Barrett and Shaun Fallat and H. Tracy Hall and Leslie Hogben and Jephian C.-H. Lin and Bryan L. Shader},
 title = {Generalizations of the Strong Arnold Property and the Minimum Number of Distinct Eigenvalues of a Graph},
 journal = {Electronic Journal of Combinatorics}, year = {2017}, volume = {24}, number = {2}, pages = {P2.40},
 doi = {10.37236/5725}
}

@article{LOS22,
 author = {Rupert H. Levene and Polona Oblak and Helena {\v{S}}migoc},
 title = {Orthogonal Symmetric Matrices and Joins of Graphs},
 journal = {Linear Algebra and its Applications}, year = {2022}, volume = {652}, pages = {213--238},
 doi = {10.1016/j.laa.2022.07.007}
}

@article{LOS25,
 author = {Rupert H. Levene and Polona Oblak and Helena {\v{S}}migoc},
 title = {Orthogonalisability of Joins of Graphs},
 journal = {Linear Algebra and its Applications}, year = {2025}, volume = {723}, pages = {162--181},
 doi = {10.1016/j.laa.2025.06.001}
}

@article{BFFK23,
 author={Wayne Barrett and Shaun Fallat and Veronika Furst and Franklin Kenter and Shahla Nasserasr and Brendan Rooney and Michael Tait and Hein van der Holst},
 title={Sparsity of graphs that allow two distinct eigenvalues},
 journal={Linear Algebra and its Applications}, volume={674}, pages={377--395}, year={2023},
 doi={10.1016/j.laa.2023.06.004}
}

@article{BFFR23,
 author={Wayne Barrett and Shaun Fallat and Veronika Furst and Shahla Nasserasr and Brendan Rooney and Michael Tait},
 title={Regular graphs of degree at most four that allow two distinct eigenvalues},
 journal={Linear Algebra and its Applications}, volume={679}, pages={127--164}, year={2023},
 doi={10.1016/j.laa.2023.09.012}
}

@article{AFK23,
 author={Aida Abiad and Shaun M. Fallat and Mark Kempton and Rupert H. Levene and Polona Oblak and Helena {\v{S}}migoc and Michael Tait and Kevin N. Vander Meulen},
 title={Bordering of symmetric matrices and an application to the minimum number of distinct eigenvalues for the join of graphs},
 journal={Linear Algebra and its Applications}, volume={679}, pages={104--126}, year={2023},
 doi={10.1016/j.laa.2023.09.013}
}

@article{LOS24,
 author={Rupert H. Levene and Polona Oblak and Helena {\v{S}}migoc},
 title={Distinct eigenvalues are realizable with generic eigenvectors},
 journal={Linear and Multilinear Algebra}, volume={72}, number={12}, pages={2054--2068}, year={2024},
 doi={10.1080/03081087.2023.2232090}
}

@article{CK25,
 author={Eric Culver and Mark Kempton},
 title={Graph products that allow two distinct eigenvalues},
 journal={Electronic Journal of Linear Algebra}, volume={41}, pages={367--379}, year={2025},
 doi={10.13001/ela.2025.9277}
}

@article{FM22Threshold,
 author={Shaun Fallat and Seyed Ahmad Mojallal},
 title={On the minimum number of distinct eigenvalues of a threshold graph},
 journal={Linear Algebra and its Applications}, volume={642}, pages={1--29}, year={2022},
 doi={10.1016/j.laa.2022.02.001}
}

@article{Allem25Threshold,
 author={L. Emilio Allem and Carlos Hoppen and Jo{\~a}o Lazzarin and Lucas Siviero Sibemberg and Fernando Colman Tura},
 title={The minimum number of distinct eigenvalues of a threshold graph is at most {4}},
 journal={Linear Algebra and its Applications}, volume={726}, pages={32--53}, year={2025},
 doi={10.1016/j.laa.2025.07.020}
}

@misc{BFP26Threshold,
 author={Jane Breen and Shaun Fallat and Johnna Parenteau},
 title={Threshold Graphs Allow Few Distinct Eigenvalues: A New Approach},
 year={2026}, note={arXiv preprint arXiv:2609.10623, September 2026},
 eprint={2609.10623}, archivePrefix={arXiv}, primaryClass={math.CO}
}

@misc{AFH26Cograph,
 author={Luiz Emilio Allem and Martin F{\"u}rer and Carlos Hoppen and Lucas Siviero Sibemberg and Vilmar Trevisan},
 title={On the minimum number of eigenvalues of matrices associated with cographs},
 year={2026}, note={arXiv preprint arXiv:2602.07282, February 2026},
 eprint={2602.07282}, archivePrefix={arXiv}, primaryClass={math.CO}
}
\end{document}